%% file: main.tex
\documentclass[12pt,draftcls,onecolumn]{IEEEtran}

\usepackage{amsmath,amsfonts, amssymb,color}
\usepackage{epstopdf}
\usepackage{nomencl}
\usepackage{enumerate}
\usepackage{times} 
\usepackage{latexsym}
\usepackage{xspace}
\usepackage{amssymb,amsmath,amsfonts,color}
\usepackage{cite}
\usepackage{caption}
\usepackage{tikz,subcaption,cite}
\usetikzlibrary{decorations.pathmorphing}

\newtheorem{theorem}{Theorem}
\newtheorem{corollary}{Corollary}
\newtheorem{proposition}{Proposition}
\newtheorem{lemma}{Lemma}
\newtheorem{remark}{Remark}
\newtheorem{definition}{Definition}
\newtheorem{example}{Example}
\DeclareMathOperator{\trace}{trace}
\DeclareMathOperator{\diag}{diag}

\definecolor{cyan}{rgb}{0.0, 1.0, 1.0}

\begin{document}
\title{\LARGE \bf  Robustness of Leader-Follower Networked Dynamical Systems}

\author{Mohammad Pirani, Ebrahim Moradi Shahrivar, Baris Fidan and Shreyas Sundaram 
\thanks{This material is based upon work supported in part by the Natural Sciences and Engineering Research Council of Canada (NSERC). M. Pirani and B. Fidan are with the Department of Mechanical and Mechatronics Engineering at the University of Waterloo, Waterloo, ON, Canada.  E-mail: \texttt{\{mpirani, fidan\}@uwaterloo.ca}. E. M. Shahrivar is with the Department of Electrical and Computer Engineering at the University of Waterloo, Waterloo, ON, Canada. E-mail: {\texttt{emoradis@uwaterloo.ca}}. S. Sundaram is with the School of Electrical and Computer Engineering at Purdue University, W. Lafayette, IN, USA.  E-mail: {\texttt{sundara2@purdue.edu}}. }
}

\maketitle

\thispagestyle{empty}
\pagestyle{empty}


\begin{abstract}

We  present a graph-theoretic approach to analyzing the robustness of  leader-follower consensus dynamics to disturbances and time delays. Robustness to disturbances is captured via the system $\mathcal{H}_2$ and $\mathcal{H}_{\infty}$ norms and robustness to time delay is defined as the maximum allowable delay for the system to remain asymptotically stable. Our analysis is built on understanding certain spectral properties of the grounded Laplacian matrix that play a key role in such dynamics. Specifically, we give graph-theoretic bounds on the extreme eigenvalues of the grounded Laplacian matrix which quantify the impact of disturbances and time-delays on the leader-follower dynamics. We then provide tight characterizations of these robustness metrics in Erdos-Renyi random graphs and random regular graphs. Finally, we view robustness to disturbances and time delay as network centrality metrics, and provide conditions under which a leader in a network optimizes each robustness objective. Furthermore, we propose a sufficient condition under which a single leader optimizes both robustness objectives simultaneously. 
\end{abstract}

\input{intro}

\input{notation}

\input{Body}

\section{Summary and Conclusions}
\label{sec:conc}

We investigated  the robustness of leader-follower consensus dynamics to uncertainty and time delay. The analysis was built on the spectrum of the grounded Laplacian matrix, which allowed us to provide tight characterizations of these robustness metrics for random graphs. Moreover, we analyzed the problem of leader selection  to optimize each robustness metric and provided a sufficient condition for a single leader to optimize all robustness metrics simultaneously. An interesting avenue for  future work is to extend the leader selection for these robustness metrics for the case of multiple leaders and analyze these robustness metrics for other classes of networks.


\bibliographystyle{IEEEtran}
\bibliography{refs}

\end{document}

%% file: intro.tex
\section{Introduction}
\label{sec:intro}

In networked dynamical systems, the influence of each individual agent on the global dynamics is determined by: (i) the location of the agent in the network, (ii) its local behavior and dynamics, and (iii) the overall network structure. One class of local behavior that has received attention in the literature (particularly in the context of social and economic networks) is the notion of agent {\it stubbornness} \cite{Friedkin2, Frasca, Ghaderi13, Ozdaglar, Acemoglu2013opinion}, where an agent influences other agents but is not affected in return.  Such agents have also been studied from the perspective of acting as {\it leaders} in multi-agent systems \cite{Rahmani,Chapman, Ni}. In this paper, our goal is to study the impact of {\it disturbances} (such as noise, faults, attacks and other external inputs) and {\it time-delays} in networked dynamical systems that contain stubborn (or leader) agents.  

From a control-theoretic viewpoint, robustness to disturbances is investigated by studying their impact on the state or output of the system, and is often quantified via system $\mathcal{H}_{2}$ or $\mathcal{H}_{\infty}$ norms. In this direction, a large  literature has recently investigated the robustness of networked dynamical systems to disturbances from an input-output standpoint \cite{Leonard,Fitch2,Bamieh2,Jovanovic,Patterson1,Clark5,Summers,Scardovi,Scardovi2,Siami,Siami,Hinf1,Hinf2,Hinf3}.
Similarly, robustness of a consensus network to time-delays in the communication between agents is quantified in terms of the maximum allowable time-delay for the system to remain asymptotically stable \cite{Olfati2}. 
 
As we describe in the paper, the robustness metrics of interest are a function of the spectrum of the grounded Laplacian matrix (obtained by removing certain rows and columns from the Laplacian matrix \cite{Barooah,Miekkala,PiraniSundaramArxiv}). 
Hence, one of the main contributions of this paper, which is an extended version of the conference paper \cite{piranicdc}, is to propose graph-theoretic bounds on the extreme eigenvalues of this matrix. These eigenvalue bounds consequently provide bounds on the robustness metrics in general graphs, and tight bounds on such metrics in random graphs.  

After characterizing graph-theoretic bounds on these robustness metrics, we turn our attention to selecting leaders in the network in order to optimize robustness.  Leader selection algorithms for multi-agent systems have attracted much attention in recent years \cite{Patterson2, Clarkbook}, and optimal leaders for certain metrics have been characterized in terms of either known network centrality measures, or via the introduction of new centrality measures \cite{Pasqualetti,Leonard,ACC}. We contribute to this literature by investigating the leader selection problem in a given network to optimize network robustness to disturbances and time delay. More specifically, we provide conditions under which a leader in a network  optimizes each robustness objective. Furthermore, we propose a sufficient graph-theoretic condition for a particular leader to optimize all of the robustness objectives simultaneously. 

The paper is organized as follows.  We start by introducing our notation in Section~\ref{sec:notation}, and in Section~\ref{sec:influence}, we formally state the leader-follower consensus dynamics that we will be studying in this paper.  There, we also describe how the spectrum of the grounded Laplacian matrix plays a role in the robustness of such dynamics to disturbances and time-delays.  In Section~\ref{sec:grounded_laplacian}, we provide graph-theoretic characterizations of the salient eigenvalues of the grounded Laplacian matrix; we then use these characterizations to provide bounds on the robustness metrics in general graphs (Section~\ref{sec:application}) and in random graphs (Section~\ref{sec:random}).  We consider the problem of selecting leaders to optimize these metrics (viewed in terms of centrality measures) in Section~\ref{sec:robustleader}, validate our analysis via simulations  in Section~\ref{sec:simulations}, and conclude in Section~\ref{sec:conc}.

%% file: notation.tex
\section{Notation}
\label{sec:notation}

We denote an undirected graph by  $\mathcal{G}=\{\mathcal{V},\mathcal{E}\}$,  where $\mathcal{V} = \{v_1, v_2, \ldots, v_n\}$ is a set of nodes (or vertices) and $\mathcal{E} \subset \mathcal{V}\times\mathcal{V}$ is the set of edges.  The neighbors of node $v_i \in \mathcal{V}$ are given by the set $\mathcal{N}_i = \{v_j \in \mathcal{V} \mid (v_i, v_j) \in \mathcal{E}\}$. The  adjacency matrix of the graph is given by a symmetric and binary $n \times n$  matrix $A$, where element $A_{ij}=1$ if $(v_i, v_j) \in \mathcal{E}$ and zero otherwise.  The degree of node $v_i$ is denoted by  $d_i \triangleq \sum_{j=1}^nA_{ij}$. For a given set of nodes $X \subset \mathcal{V}$, the {\it edge-boundary} (or just boundary) of the set is given by $\partial{X} \triangleq \{(v_i,v_j) \in \mathcal{E} \mid v_i \in X, v_j \in \mathcal{V}\setminus{X}\}$.    The Laplacian matrix of the graph is given by $L \triangleq D - A$, where $D = \diag(d_1, d_2, \ldots, d_n)$.  The eigenvalues of the Laplacian are real and nonnegative, and are denoted by $0 = \lambda_1(L) \le \lambda_2(L) \le \ldots \le \lambda_n(L)$.   For a given subset $\mathcal{S} \subset \mathcal{V}$ of nodes (which we term {\it grounded nodes}), the {\it grounded Laplacian} induced by $\mathcal{S}$ is denoted by $L_g(\mathcal{S})$ or simply $L_g$, and is obtained by removing the rows and columns of $L$ corresponding to the nodes in $\mathcal{S}$.

%% file: Body.tex
\section{Problem Statement}
\label{sec:influence}

Consider a connected network  consisting of $n$ agents $\mathcal{V} = \{v_1, v_2,\ldots, v_n\}$. The set of agents is partitioned into a set of followers $\mathcal{F}$, and a set of leaders\footnote{These agents may also be referred to as {\it{anchors}} \cite{Rahmani} or {\it{stubborn agents}}\cite{Ghaderi13},\cite{piranicdc} depending on the context.} $\mathcal{S}$. Each agent $v_i$ has a scalar and real valued state $x_i(t)$, where $t$ is the time index. The state of each follower agent $v_j\in \mathcal{F}$ evolves based on the interactions  with its neighbors as
\begin{align}
\dot{x}_{j}(t)&=\sum_{v_i\in \mathcal{N}_j}(x_i(t)-x_j(t)).
\label{eqn:partial}
\end{align}
The state of the leaders (which should be tracked by the followers)  is assumed to be constant\footnote{The results in this paper can be  extended to the case where the state of the leaders are time-varying  \cite{Rahmani}, and  given by $\dot{x}_{\mathcal{S}}(t) = u(t)$.} and thus
\begin{equation}
\dot{x}_{j}(t) = 0, \enspace \forall v_j \in \mathcal{S}.
\label{eqn:fully}
\end{equation}

If the graph is connected, the  states of the follower agents will converge to  some convex combination of the states of the leaders \cite{Clark1}. We assume without loss of generality that the leader agents are placed last in the ordering of the agents. Aggregating the states of all followers  into a vector $x_\mathcal{F}(t) \in \mathbb{R}^{n-|\mathcal{S}|}$, and the states of all  leaders into a vector $x_{\mathcal{S}}(t)\in \mathbb{R}^{|\mathcal{S}|}$ (note that $x_{\mathcal{S}}(t) = x_{\mathcal{S}}(0)$ for all $t \ge 0$), equations \eqref{eqn:partial} and \eqref{eqn:fully} yield the following dynamics 
\begin{equation}
\begin{bmatrix}
      \dot{x}_\mathcal{F}(t)          \\[0.3em]
       \dot{x}_{\mathcal{S}}(t) 
     \end{bmatrix}=-\underbrace{\begin{bmatrix}
       L_g & L_{12}           \\[0.3em]
       L_{21} & L_{22}           
     \end{bmatrix}}_L\begin{bmatrix}
      {x}_F(t)          \\[0.3em]
       {x}_{\mathcal{S}}(t) 
     \end{bmatrix}.
\label{eqn:mat}
\end{equation}
Given equation \eqref{eqn:fully}, we have that $ L_{21}=0$ and $ L_{22}=0 $. Hence the dynamics of the follower agents are given by 
\begin{align}
\dot{x}_\mathcal{F}(t) &= -{L}_gx_\mathcal{F}(t) + L_{12}x_{\mathcal{S}}(0).
\label{eqn:partial4}
\end{align}
 Here, $L_g$ is the grounded Laplacian induced by the leaders, representing the interaction between the followers. The submatrix $L_{12}$ of the graph Laplacian captures the influence of the leaders on the followers.

 \begin{remark}
For the case where the underlying network is connected and there exists at least one leader, the matrix ${L}_g$ is a diagonally dominant matrix with at least one strictly diagonally dominant row. Hence, from \cite{Horn} it is a positive definite matrix and  in this case,  the dynamics given by \eqref{eqn:partial4} will be asymptotically stable and the  convergence rate is determined by the  smallest eigenvalue of ${L}_g$.  Moreover ${L}_g^{-1}$ is a nonnegative matrix and based on the Perron-Frobenius theorem \cite{Horn} its largest eigenvalue $\lambda_{n-|\mathcal{S}|}({L}_g^{-1})$ has an eigenvector $\mathbf{x}$ with nonnegative components.  Thus, the smallest eigenvalue $\lambda_1(L_g)$ of $L_g$ also has an eigenvector $\mathbf{x}$ with nonnegative components.
\label{rem:deffin}
\end{remark}

 In this paper, we will consider the impact of two perturbations to the above nominal dynamics\footnote{We analyze these two cases separately, since the objective of this study is to show the explicit dependency of these two robustness metrics on the spectrum of $L_g$.}: 
 
 \begin{itemize}
 \item The  case where  the update rule of each follower agent $v_j \in \mathcal{F}$  is affected by a disturbance $w_j(t)$. In this case, we extend \eqref{eqn:partial4} to  
\begin{align}
\dot{{x}}_\mathcal{F}(t) &= -{L}_g{x}_\mathcal{F}(t)+ L_{12}x_{\mathcal{S}}(0) + w(t).
\label{eqn:padgrtial43}
\end{align}
 Here $w(t)$ is a vector representing the disturbances. It is assumed that the leader agents are unaffected by the disturbances. This is a reasonable assumption due to the fact that they do not update their state.
 
 \item  The case where the communication  between the agents is affected by some time delay.  In this case,  we have the  dynamics
 \begin{align}
\dot{{x}}_\mathcal{F}(t) &= -{L}_g{x}_\mathcal{F}(t-\tau) + L_{12}x_{\mathcal{S}}(0),
\label{eqn:partiweal43}
\end{align}
 where $0<\tau \leq \tau_{max}$ for  some $\tau_{max}>0$. 
 \end{itemize}
 
\begin{remark}
If each agent has instantaneous access to its own state, the dynamics have the form
\begin{equation}
\dot{{x}}_\mathcal{F}(t) = -D_g{x}_\mathcal{F}(t)+ A_g{x}_\mathcal{F}(t-\tau) + L_{12}x_{\mathcal{S}}(0),
\label{eqn:offdiag}
\end{equation}
where $L_g=D_g-A_g$. In this case  since all of the principal minors of $L_g$ are nonnegative and $L_g$ is non-singular, \eqref{eqn:offdiag} is asymptotically stable independent of the magnitude of the delays in the off-diagonal terms of $L_g$ (Theorem 1 in \cite{Hofbauer}).
\end{remark}

 In the following subsections, we analyze the robustness of system \eqref{eqn:padgrtial43} to the disturbances and \eqref{eqn:partiweal43} to the time delay.

\subsection{Robustness of \eqref{eqn:padgrtial43} to Disturbances}
\label{sec:robust_disturbances}

 Let $\bar{x}_\mathcal{F}(t)$ be the state of \eqref{eqn:padgrtial43} when $w(t)=0$, and define the error between the nominal and disturbed state as  $e(t)={x}_\mathcal{F}(t)-\bar{x}_\mathcal{F}(t)$. The transfer function from the disturbance $w(t)$ to $e(t)$ is obtained from \eqref{eqn:padgrtial43} as $G(s)= (sI+L_g)^{-1}$. In order to discuss the robustness of \eqref{eqn:padgrtial43} to disturbances, a typical approach (e.g. \cite{Bamieh,Jovanovic,piranicdc}) is to consider system $\mathcal{H}_2$ and $\mathcal{H}_{\infty}$ norms, defined as \cite{Doyle}
\begin{align}
||G||_2 &\triangleq \left( \frac{1}{2\pi}\trace\int_0^{\infty}G^*(j\omega)G(j\omega)d\omega \right)^{\frac{1}{2}},\nonumber \\
||G||_{\infty} &\triangleq \sup_{\omega\in \mathbb{R}}{\lambda_{max}^{\frac{1}{2}}(G^*(j\omega)G(j\omega))}.
\end{align}

The system $\mathcal{H}_2$ norm can also be calculated based on the controllability Gramian $\mathcal{W}_c$, which is the solution of a Lyapunov equation.  In particular, for  the error dynamics of \eqref{eqn:padgrtial43} we have $||G||_2^2=\trace \mathcal{W}_c$, which becomes \cite{Bamieh2}
\begin{equation}
||G||_2=\left(\frac{1}{2}\trace(L_g^{-1})\right)^{\frac{1}{2}}. 
\label{eqn:htwo}
\end{equation}

For the system $\mathcal{H}_{\infty}$ norm of the error dynamics of \eqref{eqn:padgrtial43}, we present the following proposition.

\begin{proposition}
The system  $\mathcal{H}_{\infty}$ norm  of the error dynamics of \eqref{eqn:padgrtial43} is 
\begin{equation}
||G||_{\infty}=\frac{1}{\lambda_1({L}_g)}.
\label{eqn:hinfti}
\end{equation}
\label{prop:hinf}
\end{proposition}

\begin{IEEEproof}
We have $G(j\omega)=(j\omega I+{L}_g)^{-1}$, which gives
\begin{align}
G^*(j\omega)G(j\omega)&=(-j\omega I+{L}_g)^{-1}(j\omega I+{L}_g)^{-1} \nonumber \\
&=\underbrace{(\omega^2I+{L}_g^2)^{-1}}_{\mathcal{C}^{-1}}.
\end{align}
We know that $\mathcal{C}>0$ (positive definite) which yields $\mathcal{C}^{-1}>0$. Thus finding  $\sup_{\omega} \lambda_{n-|\mathcal{S}|}(\mathcal{C}^{-1})$  is equivalent to finding  $\inf_{\omega} \lambda_{1}(\mathcal{C})$. Since $\lambda_{1}(\mathcal{C})=\omega^2+\lambda_{1}({L}_g^2)$, we have $\inf_{\omega} \lambda_{1}(\mathcal{C})=\lambda_{1}({L}_g^2)$, proving the proposition.  
\end{IEEEproof}

Based on Proposition \ref{prop:hinf} and Remark \ref{rem:deffin}, minimizing the $\mathcal{H}_{\infty}$ norm of the error dynamics of \eqref{eqn:padgrtial43} is equivalent to maximizing the convergence rate of \eqref{eqn:padgrtial43}.  This equivalence between the two metrics will be revisited in Section \ref{sec:robustleader}. 

\begin{remark}
The recent literature  mainly works with $\frac{1}{2}\trace(L_g^{-1})$ instead of its square root, and refers to this metric as the {\it {network  disorder}} \cite{Bamieh2}. To maintain consistency, we also adopt this terminology and refer to $\frac{1}{2}\trace(L_g^{-1})$ as $\mathcal{H}_2$ disorder and  $\frac{1}{\lambda_1({L}_g)}$ as $\mathcal{H}_{\infty}$ disorder.
\end{remark}

\subsection{Robustness of \eqref{eqn:partiweal43} to Time Delay}

The other robustness metric that we analyze in this paper is the robustness of \eqref{eqn:partiweal43} to  time delay. The following theorem gives a necessary and sufficient condition for the asymptotic stability of \eqref{eqn:partiweal43}. 

\begin{theorem}[\cite{Buslowicz}]
The dynamical system  \eqref{eqn:partiweal43} is asymptotically stable if and only if
\begin{equation}\label{thm:ineq1}
\tau_{max}< \frac{\tan^{-1}\left(\frac{Re(\lambda_i(L_g))}{Im(\lambda_i(L_g))}\right)}{\lambda_i(L_g)},
\end{equation}
for all $i=1,2, ..., n-|\mathcal{S}|$. 
\label{thm:delll}
\end{theorem}

Since the eigenvalues are real, the numerator of the right hand side term in   \eqref{thm:ineq1} is  $\frac{\pi}{2}$, and thus a necessary and sufficient condition for asymptotic stability of \eqref{eqn:partiweal43} is 
\begin{equation}
\tau_{max}<\frac{\pi}{2\lambda_{max}(L_g)}=\frac{\pi}{2\lambda_{n-|S|}(L_g)}.
\label{eqn:ti}
\end{equation}

In this paper, we refer to the quantity $\hat{\tau}_{max}=\frac{\pi}{2\lambda_{n-|S|}(L_g)}$ as the {\it{delay threshold}}.

The above characterization of the robustness of \eqref{eqn:padgrtial43} to disturbances, based on system $\mathcal{H}_2$ and $\mathcal{H}_{\infty}$ norms in \eqref{eqn:htwo} and \eqref{eqn:hinfti}, and the robustness of  \eqref{eqn:partiweal43} to delay given by the quantity in \eqref{eqn:ti}, illustrates the role that the spectrum of the grounded Laplacian matrix $L_g$ plays in such robustness metrics. Thus, we analyze the spectrum of the grounded Laplacian matrix in this paper, and consequently give bounds on the above robustness metrics.  More specifically, our contributions are as follows.

\begin{itemize}
\item We extend existing  bounds on the smallest eigenvalue, $\lambda_1$, and provide new bounds on the largest eigenvalue, $\lambda_{n-|\mathcal{S}|}$, of the grounded Laplacian matrix.

\item Based on the graph-theoretic bounds on the extreme eigenvalues of $L_g$, we present graph-theoretic necessary and sufficient conditions for robustness of leader-follower dynamics to disturbances and time delay. 

\item We characterize the system $\mathcal{H}_2$ disorder, $\mathcal{H}_{\infty}$ disorder, and the robustness to delay ($\hat{\tau}_{max}$) in random graphs. 

\item We look at these robustness metrics for the disturbance and time delay as different network centrality metrics  and give sufficient conditions for a node in a network to be the best leader in the sense of minimizing $\mathcal{H}_2$ disorder, minimizing $\mathcal{H}_{\infty}$ disorder (or maximizing convergence rate) and maximizing $\hat{\tau}_{max}$ simultaneously. 
\end{itemize}


\section{On the Spectrum of the Grounded Laplacian Matrix}
\label{sec:grounded_laplacian}

In this section, we present graph-theoretic bounds on  the smallest eigenvalue and the spectral radius (largest eigenvalue) of the grounded Laplacian matrix.  

\subsection{Smallest Eigenvalue of $L_g$}

There is a vast literature  dedicated to analyzing the spectrum  of the Laplacian matrix \cite{Mohar}, \cite{Anderson}, \cite{Chung}. The following theorem gives  bounds on $\lambda_1(L_g)$ (the smallest eigenvalue of the grounded Laplacian matrix) based on  graph theoretic properties. 

\begin{theorem}
Consider a connected graph $\mathcal{G}=\{\mathcal{V},\mathcal{E}\}$ with a set of leaders $\mathcal{S} \subset \mathcal{V}$.  Let $L_g$ be the grounded Laplacian matrix induced by $\mathcal{S}$, and for each $v_i \in \mathcal{F}$, let $\beta_i$ be the number of leaders in follower $v_i$'s neighborhood. Then
\begin{multline}
\max\left\{ \min_{i\in \mathcal{V}\setminus \mathcal{S}}\beta_i, \left(\frac{|\partial \mathcal{S}|}{n-|\mathcal{S}|}\right)x_{min}\right\} \leq \lambda_1({L}_g) \leq 
\min_{\emptyset \ne X\subseteq \mathcal{V\setminus \mathcal{S}}} \frac{|\partial X|}{|X|} \leq \frac{|\partial \mathcal{S}|}{n-|\mathcal{S}|}\leq \max_{i\in \mathcal{V}\setminus \mathcal{S}}\beta_i,
\label{eqn:maineq}
\end{multline}
where $x_{min}$ is the smallest eigenvector component of $\mathbf{x}$, a nonnegative eigenvector corresponding to $\lambda_1({L}_g)$.\footnote{Throughout the paper, we take eigenvector $\mathbf{x}$ to be normalized such that its largest component is $x_{max}=1$.}
\label{thm:main}
\end{theorem}

\begin{IEEEproof}
The lower bound $\frac{|\partial \mathcal{S}|}{n-|\mathcal{S}|}x_{min}$ and the tightest and the second tightest upper bounds are given in Theorem 1 in \cite{PiraniSundaramArxiv}.
 The extreme upper bound is due to the fact that $\sum_{i=1}^{n-|\mathcal{S}|}\beta_i=|\partial \mathcal{S}|$ which gives $\frac{|\partial \mathcal{S}|}{n-|\mathcal{S}|}\leq \max_{i\in \mathcal{V}\setminus \mathcal{S}}\beta_i$.
 For the lower bound $\min_{i\in \mathcal{V}\setminus \mathcal{S}}\beta_i$, we left-multiply the eigenvector equation $\lambda_1\textbf{x}={L}_g\textbf{x}$ by the vector consisting of all 1's, and use the fact that $\mathbf{1}^TL_g=[\beta_1, \beta_2, ..., \beta_{n-|\mathcal{S}|}]$ to get
\begin{equation*}
  \lambda_1({L}_g)\sum_{i=1}^{n-|\mathcal{S}|}x_i=\sum_{i=1}^{n-|\mathcal{S}|}\beta_ix_i\geq \min_{i\in \mathcal{V}\setminus \mathcal{S}}\beta_i\sum_{i=1}^{n-|\mathcal{S}|}x_i,
 \end{equation*}
since $\mathbf{x}$ is nonnegative, which gives $\lambda_1({L}_g)\geq \min_{i\in \mathcal{V}\setminus \mathcal{S}}\beta_i$ as required.
\end{IEEEproof}

The following lemma from \cite{PiraniSundaramArxiv}  provides a sufficient condition under which the smallest component of the eigenvector corresponding to $\lambda_1({L}_g)$ goes to $1$ and consequently the bound \eqref{eqn:maineq} becomes tight. 

\begin{lemma}[\cite{PiraniSundaramArxiv}]
Let $\mathbf{x}$ be a nonnegative eigenvector corresponding to the smallest eigenvalue of ${L}_g$.  Then the smallest eigenvector component of $\mathbf{x}$ satisfies
\begin{equation}
x_{min} \ge 1 - \frac{2\sqrt{|\mathcal{S}||\partial \mathcal{S}|)}}{\lambda_2({\bar{L}})},
\label{eqn:eige12}
\end{equation}
where $\bar{L} \in \mathbb{R}^{(n-|\mathcal{S}|)\times (n-|\mathcal{S}|)}$ is the Laplacian matrix formed by removing the leaders and their incident edges.
 \label{lem:mainlem}
\end{lemma}

Thus, in networks where the number of leaders (and their incident edges)  grows slowly compared to the algebraic connectivity of the network induced by the followers, the bounds on the smallest eigenvalue in \eqref{eqn:maineq} become tight. 

In order to obtain another condition under which $\lambda_1({L}_g)$ is bounded away from zero, we use the following definition. 

\begin{definition}[\cite{Kutten}]
A subset of vertices $\mathcal{X}\subset \mathcal{V}$ is an $f$- dominating set if each vertex $v_i\in \mathcal{V}\setminus \mathcal{X}$ is connected to at least $f$ vertices in $\mathcal{X}$.
\end{definition}

Based on the above definition and the lower bound $\min_{i\in \mathcal{V}\setminus \mathcal{S}}\{\beta_i\}$ in \eqref{eqn:maineq}, we have the following corollary.

\begin{corollary}
If the set of leaders is an $f$-dominating set, then  $\lambda_1({L}_g)\geq f$, regardless of the  connectivity of the network.
\label{cor:akjdbvao}
\end{corollary}

The following proposition introduces a condition under which $\lambda_1({L}_g)$ remains unchanged when some edges are added or removed.

\begin{proposition}
Consider a connected graph $\mathcal{G}=\{\mathcal{V},\mathcal{E}\}$ with a set of leaders $\mathcal{S} \subset \mathcal{V}$.  Let $L_g$ be the grounded Laplacian matrix induced by $\mathcal{S}$. If each follower is connected to exactly $\beta$ leaders, then regardless of the interconnection topology inside $\mathcal{F}$ or $\mathcal{S}$, we have $\lambda_1(L_g)=\beta$. Moreover in this case $\lambda_1(L_g)$ strictly decreases when any edge  $(v_i,v_j)\in \mathcal{E}$ is removed, where $v_i\in \mathcal{F}$ and $v_j\in \mathcal{S}$. 
\label{prop:robus}
\end{proposition}

\begin{IEEEproof}
The proof of the first part is clear due to \eqref{eqn:maineq} since $\min_{i\in \mathcal{F}}\{\beta_i\}=\max_{i\in \mathcal{F}}\{\beta_i\}=\beta$. Furthermore, in this case if we remove an edge between $\mathcal{F}$ and $\mathcal{S}$, then based on \eqref{eqn:maineq} we have $\lambda_1(L_g) \leq \frac{|\partial \mathcal{S}|}{n-|\mathcal{S}|}=\frac{(n-|\mathcal{S}|-1)\beta+\beta-1}{n-|\mathcal{S}|}<\beta$, which proves the claim.
\end{IEEEproof}

Proposition \ref{prop:robus}  is important from two aspects. First, it gives freedom in designing connections between the follower agents. Second, it introduces a notion of robustness of the network under edge failures within the set of follower (or leader) agents.

\subsection{Spectral Radius of $L_g$}

Bounds on the spectral radius of the Laplacian matrix are discussed in \cite{Anderson, Shiiii}. Here we discuss  graph theoretic bounds on the spectral radius,  $\lambda_{n-|\mathcal{S}|}(L_g)$, of the grounded Laplacian matrix. We start with the definition of the {\it incidence matrix}.
 
 \begin{definition}
 Given a connected  graph $\mathcal{G}=\{\mathcal{V},\mathcal{E}\}$, an orientation of the graph $\mathcal{G}$ is defined by assigning  a direction (arbitrarily) to each edge in $\mathcal{E}$. For graph $\mathcal{G}$ with $m$ edges, numbered as $e_1, e_2, ..., e_m$, its node-edge incidence matrix $\mathcal{B}(\mathcal{G})\in \mathbb{R}^{n\times m}$ is defined as
$$[\mathcal{B}(\mathcal{G})]_{kl}=
  \begin{cases}
    1       & \quad  \text{if node $k$ is the head of edge $l$},\\
   -1  & \quad \text{if node $k$ is the tail of edge $l$},\\
   0  & \quad \text{otherwise}.\\
  \end{cases}
  $$
The graph Laplacian satisfies $L=\mathcal{B}(\mathcal{G})\mathcal{B}(\mathcal{G})^T$. 
  \end{definition}
  
Partitioning the rows of the incidence matrix into the sets of followers and leaders yields
\begin{equation}
\mathcal{B}(\mathcal{G})=\begin{bmatrix} \mathcal{B}_{\mathcal{F}} \\[0.3em] 
\mathcal{B}_{\mathcal{S}} \end{bmatrix},
\end{equation}
where $\mathcal{B}_{\mathcal{F}}\in \mathbb{R}^{(n-|\mathcal{S}|)\times m}$ and $\mathcal{B}_{\mathcal{S}}\in \mathbb{R}^{|\mathcal{S}|\times m}$. 
As a result, $L_g=\mathcal{B}_{\mathcal{F}} \mathcal{B}_{\mathcal{F}} ^T$. Defining the matrix $N=\mathcal{B}_{\mathcal{F}}^T \mathcal{B}_{\mathcal{F}}$, we have the following lemma. 
\begin{lemma}
If $\lambda_k$ is an eigenvalue of  $L_g$, it is also an eigenvalue of $N$.
\label{lem:fer}
\end{lemma}

\begin{IEEEproof}
If we have $L_g\mathbf{x}_k=\lambda_k\mathbf{x}_k$ for any  eigenvalue $\lambda_k(L_g)$ and corresponding eigenvector $\mathbf{x}_k$, then $N\mathcal{B}_{\mathcal{F}}^T\mathbf{x}_k=\mathcal{B}_{\mathcal{F}}^TL_g\mathbf{x}_k=\lambda_k\mathcal{B}_{\mathcal{F}}^T\mathbf{x}_k$.   If $\mathcal{B}_{\mathcal{F}}^T\mathbf{x}_k=\mathbf{0}$ then $L_g\mathbf{x}_k=\mathbf{0}$ which is impossible since $L_g$ is positive definite by Remark \ref{rem:deffin}. Thus $\lambda_k$ is also an eigenvalue of $N$ with eigenvector $\mathcal{B}_{\mathcal{F}}^T\mathbf{x}_k$.
\end{IEEEproof}

This leads to the following bounds on $\lambda_{n-|\mathcal{S}|}(L_g)$. 

\begin{theorem}
Consider a connected graph $\mathcal{G}=\{\mathcal{V},\mathcal{E}\}$ with a set of leaders $\mathcal{S} \subset \mathcal{V}$.  Let $L_g$ be the grounded Laplacian matrix induced by $\mathcal{S}$. The spectral radius of $L_g$ satisfies
\begin{align}
d_{max}^{\mathcal{F}} \leq \lambda_{n-|\mathcal{S}|}(L_g) 
\leq \max\left\{d_{max}^{\mathcal{F}}, \max_{\substack{(u,v)\in \mathcal{E}\\{u,v\in \mathcal{F}}}}\{d_u+d_v\} \right\},
\label{eqn:up}  
\end{align}
where $d_{max}^{\mathcal{F}}$ is the maximum degree over the follower agents.
\label{thm:bounddmax}
\end{theorem}

\begin{IEEEproof}
For the lower bound, the Rayleigh quotient inequality  \cite{Horn} indicates
$$
\lambda_{n-|\mathcal{S}|}(L_g) \geq z^T L_g z ,
$$
for all $z\in \mathbb{R}^{n-|S|}$ with $z^Tz=1$. By choosing $z=e_i$, where $e_i$ is a vector of zeros except for a single 1 at an element corresponding to a vertex with maximum degree, the lower bound is obtained. 

In order to show the upper bound, we use the property of matrix $N$ mentioned in Lemma   \ref{lem:fer}. Thus we show the same upper bound for the spectral radius of  $N$, $\lambda_{max}(N)$. We have
\begin{equation}
\lambda_{max}(N)\leq \lambda_{max}(|N|) \leq \max_{i=1, \dots, n-|S|}[|N|]_i,
\label{eqn:pf}
\end{equation}
where $[|N|]_i$ is the row sum of the $i$-th row of $|N|$. The first inequality in \eqref{eqn:pf} is due to the properties of  nonnegative matrices (Theorem 8.1.18 in \cite{Horn}) and  the second inequality  is due to the Perron-Frobenius theorem. We know that the $i$-th row of  $N$   belongs to edge $e_i$. This edge is  either connecting two follower agents or a follower with a leader. Moreover, when $e$ is an edge in $\mathcal{G}$ connecting $\{u,v\}\in \mathcal{F}$, the row sum of $|N|$ for the row corresponding to $e$ is $d_u+d_v$ \cite{Anderson}. Furthermore, when $e$ is an edge in $\mathcal{G}$ which connects $u\in \mathcal{F}$ to $v\in \mathcal{S}$, then the  row sum is $d_u$. This gives the upper bound in \eqref{eqn:up}.
\end{IEEEproof}

\begin{remark}
Unlike the case where there is no leader in the network (i.e., traditional consensus dynamics) where we always have $d_{max}< \max_{(u,v)\in \mathcal{E}}\{d_u+d_v\}$, for the upper bound in \eqref{eqn:up}, it is possible to have 
\begin{equation}
d_{max}^{\mathcal{F}} > \max_{\substack{(u,v)\in \mathcal{E}\\{u,v\in \mathcal{F}}}}\{d_u+d_v\},
\label{eqn:wrnh}
\end{equation}
and \eqref{eqn:up} becomes tight, i.e., $\lambda_{n-|\mathcal{S}|}(L_g)=d_{max}^{\mathcal{F}}$. This corresponds to the case where the vertex with maximum degree among the followers is not connected to any follower (its neighbors are all leaders) and  the degrees of the rest of the followers are small enough such that \eqref{eqn:wrnh} is satisfied. 
\label{rem:aetng}
\end{remark}

Based on Theorem \ref{thm:bounddmax} and Remark \ref{rem:aetng}, we  introduce  a class of graphs in which the bound \eqref{eqn:up} is tight, i.e., $\lambda_{n-|\mathcal{S}|}(L_g)=d_{max}^{\mathcal{F}}$. 

\begin{definition}[\cite{Godsil}]
An independent vertex set of a graph $\mathcal{G}$ is a subset of the vertices such that no two vertices in the subset are connected to each other via an edge.
\end{definition}

\begin{corollary}
If the set of followers $\mathcal{F}$ is an independent set, then  $\lambda_{n-|S|}(L_g)=d_{max}^{\mathcal{F}}$. 
\label{cor:remm}
\end{corollary}

\begin{IEEEproof}
We can prove this statement in two ways. The first is based on the proof of Theorem \ref{thm:bounddmax}: since there is no row in $N$ which belongs to an edge connecting two followers, we have $\lambda_{n-|S|}(L_g)=d_{max}^{\mathcal{F}}$. The second proof is to note that in the case where $\mathcal{F}$ is an independent set, the grounded Laplacian matrix will be a diagonal matrix and $\lambda_{n-|S|}(L_g)$ will be the largest diagonal element, namely $d_{max}^{\mathcal{F}}$.
\end{IEEEproof}

A simple example that satisfies the condition in Corollary \ref{cor:remm} is a bipartite graph in which one partition consists of the leaders and the other set contains the followers.

\section{Application of the Spectrum of $L_g$ to Network Robustness to Disturbances and Time Delay}
\label{sec:application}

In the previous section, we analyzed the spectral properties of the grounded Laplacian matrix $L_g$. In this section, we use the results from the previous section to give bounds on the network robustness metrics we identified earlier.

\subsection{Robustness to Disturbances}

By Proposition \ref{prop:hinf}, we know that the system $\mathcal{H}_{\infty}$ norm of the error dynamics of \eqref{eqn:padgrtial43} is equal to $\frac{1}{\lambda_1(L_g)}$. Hence, based on Theorem \ref{thm:main}, we have the following  bounds  for $\mathcal{H}_{\infty}$ disorder.
\begin{equation}
 \frac{1}{\max_{i\in \mathcal{F}}\{\beta_i\}} \leq \frac{n-|\mathcal{S}|}{|\partial \mathcal{S}|} \leq \frac{1}{\lambda_1(L_g)} \leq \frac{1}{\min_{i\in \mathcal{F}}\{\beta_i\}}.
\label{eqn:hinfff}
\end{equation}
Note that the upper bound is taken to be $\infty$ if $\min_{i\in \mathcal{F}}\{\beta_i\}=0$. Based on \eqref{eqn:hinfff}, for a leader-follower multi-agent system with leader set $\mathcal{S}$ and follower set $\mathcal{F}$,  a necessary  condition to have  $||G||_{\infty}\leq \gamma$  is to have $\frac{1}{\max_{i\in \mathcal{F}}\{\beta_i\}} \leq \gamma$ and a sufficient condition is to have $\frac{1}{\min_{i\in \mathcal{F}}\{\beta_i\}}\leq \gamma$. Based on Corollary \ref{cor:akjdbvao}, a sufficient  condition for $||G||_{\infty}\leq \gamma$ is that  the leader set is a $\lceil \frac{1}{\gamma} \rceil$-dominating set.

\subsection{Robustness to Time Delay}

 Based on the bounds discussed in Theorem \ref{thm:bounddmax} and \eqref{eqn:ti}, a necessary  condition for asymptotic stability of \eqref{eqn:partiweal43} is 
\begin{align}
\tau_{max} < \frac{\pi}{2 d_{max}^{\mathcal{F}}}.
\label{eqn:de}
\end{align}
Moreover, a sufficient condition for asymptotic stability of \eqref{eqn:partiweal43} is
\begin{equation}
\tau_{max}  < \frac{\pi}{2 \mathcal{M}},
\label{eqn:dwede}
\end{equation}
where $\mathcal{M}=\max\left\{d_{max}^{\mathcal{F}}, \max_{\substack{(u,v)\in \mathcal{E}\\{u,v\in \mathcal{F}}}} d_u+d_v \right\}$.

\begin{remark}
Since $ d_{max}^{\mathcal{F}} \geq 1$ for connected graphs, condition \eqref{eqn:de} implies that the delay must necessarily be strictly less than $\frac{\pi}{2}$ for \eqref{eqn:partiweal43} to be stable. 
\end{remark}

\begin{remark}
Based on Corollary \ref{cor:remm}, if the set of followers is an independent set, then  necessary and sufficient conditions \eqref{eqn:de} and \eqref{eqn:dwede} coincide, i.e., $\hat{\tau}_{max} = \frac{\pi}{2 d_{max}^{\mathcal{F}}}$.
\end{remark}

\subsection{Trade-off Between Minimizing $\mathcal{H}_{\infty}$-Norm and Maximizing $\hat{\tau}_{max}$ }

There is a trade-off between maximizing $\lambda_1(L_g)$ (minimizing $\mathcal{H}_{\infty}$ disorder) and minimizing $\lambda_{n-{\mathcal{S}}}(L_g)$ (maximizing delay threshold). The same situation is discussed for the algebraic connectivity ($\lambda_2(L)$) and the spectral radius of the Laplacian matrix in \cite{Olfati2}. In this subsection, we address this problem. More formally, we introduce the following network design problem. Suppose we are given a set of agents consisting of at least $\beta$ leaders. The objective is to design a network to solve the optimization problem,

\begin{equation}
\begin{aligned}
& \underset{X}{\text{minimize}}
& & J(X)=\lambda_{n-|\mathcal{S}|}(L_g)  \\
& \text{subject to}
& & \lambda_1(L_g) \geq \beta, \\
&&& X_i\in \{0,1\}, \\
\end{aligned}
\label{eqn:ebf}
\end{equation}
where $X_{\binom{n}{2}\times 1}$ is an indicator vector in which each edge $e_i$ in the graph is assigned to an element in $X$, i.e.,  $X_i=1$ if $e_i$ exists and $X_i=0$ otherwise. Despite the Boolean constraint $X_i\in \{0,1\}$, Proposition \ref{prop:robus} provides an efficient way to solve this problem, namely  connecting each follower to exactly $\beta$ leaders.\footnote{ Note that based on Theorem \ref{thm:main}, $\lambda_1(L_g)$ is upper bounded by the total number of leaders, and thus at least $\beta$ leaders are needed to make the problem feasible.} Based on that proposition, in this case the value of $\lambda_1(L_g)$ will be independent of the interconnections between the followers. If we make the set of follower agents an independent set,  from Corollary \ref{cor:remm}, we have $\lambda_{n-|{\mathcal{S}}|}(L_g)=d_{max}^{\mathcal{F}}=\beta$. In other words, this design makes $L_g$ a diagonal matrix whose diagonal elements are $\beta$, and the optimal solution is attained.

\section{Robustness in Random Graphs}
\label{sec:random}

In this section we discuss  $\mathcal{H}_2$ and  $\mathcal{H}_{\infty}$ disorders and the delay threshold $\hat{\tau}_{max}$ when the underlying network structure is a random graph. We analyze two well known random graphs, namely Erdos-Renyi (ER) random graphs and random regular graphs (RRG).

\subsection{Erdos-Renyi Random Graphs}

\begin{definition}
An Erdos-Renyi (ER) random graph $\mathcal{G}(n,p)$ is a graph on $n$ nodes, where each edge between two distinct nodes is present independently with probability $p$ (which could be a function of $n$).  We say that a graph property holds {\it asymptotically almost surely} if the probability of drawing a graph with that property goes to $1$ as $n \rightarrow \infty$.  Let $\Omega_n$ be the set of all undirected graphs on $n$ nodes. For a given graph function $f: \Omega_n \rightarrow \mathbb{R}_{\ge 0}$ and another function $g: \mathbb{N} \rightarrow \mathbb{R}_{\ge 0}$, we say $f(\mathcal{G}(n,p)) \le (1+o(1))g(n)$ asymptotically almost surely if there exists some function $h(n) \in o(1)$ such that $f(\mathcal{G}(n,p)) \le (1+h(n))g(n)$ with probability tending to $1$ as $n \rightarrow \infty$. Lower bounds of the above form have an essentially identical definition.
\end{definition}

\subsubsection{Network Disorder in ER Random Graphs}
Before discussing  network disorder in ER random graphs we recall the following theorem for the smallest eigenvalue of the grounded Laplacian in such graphs.

\begin{theorem}[\cite{PiraniSundaramArxiv}]
Consider the Erdos-Renyi random graph $\mathcal{G}(n,p)$, where the edge probability $p$ satisfies
$p(n) \ge \frac{c\ln{n}}{n}$, for constant $c>1$.  Let $\mathcal{S}$ be a set of grounded nodes chosen uniformly at random with $|\mathcal{S}| = o(\sqrt{np})$. Then the smallest eigenvalue $\lambda_1(L_g)$ of the grounded Laplacian satisfies $(1-o(1))|\mathcal{S}|p \le \lambda_1(L_g) \le (1+o(1))|\mathcal{S}|p$ asymptotically almost surely.
\label{thm:erdos1}
\end{theorem}

The above theorem covers a broad range of edge-probability functions, and includes constant $p$ as a special case. Based on the above theorem, we obtain the following result for $\mathcal{H}_{\infty}$ disorder  in random graphs.

\begin{theorem} 
Consider a random graph $\mathcal{G}(n,p)$ with  $p(n) \ge \frac{c\ln{n}}{n}$, for constant $c>1$.  Let $\mathcal{S} \subset \mathcal{V}$ be a set of grounded nodes chosen uniformly at random with $|\mathcal{S}| = o(\sqrt{np})$. Then for $\mathcal{H}_{\infty}$ disorder we have 
\begin{equation}
(1 - o(1))\frac{1}{|\mathcal{S}|p} \le \frac{1}{\lambda_1(L_g)} \le(1 + o(1))\frac{1}{|\mathcal{S}|p},
\label{eqn:per4}
\end{equation}
asymptotically almost surely. 
\label{thm:cohh}
\end{theorem}

Moreover, we have the following result for $\mathcal{H}_{2}$ disorder  in random graphs  with constant edge probability $p$.

\begin{theorem} 
Consider a random graph $\mathcal{G}(n,p)$ with constant $p$.  Let $\mathcal{S} \subset \mathcal{V}$ be a set of grounded nodes chosen uniformly at random with $|\mathcal{S}|=o(\sqrt{n})$. Then for $\mathcal{H}_{2}$ disorder we have 
\begin{equation}
(1 - o(1))\frac{|\mathcal{S}|+1}{2|\mathcal{S}|p} \le \frac{1}{2}\trace(L_g^{-1}) \le(1 + o(1))\frac{|\mathcal{S}|+1}{2|\mathcal{S}|p},
\label{eqn:per44}
\end{equation}
asymptotically almost surely. 
\label{thm:coh}
\end{theorem}

\begin{IEEEproof}
For each node $v_i \in \mathcal{V}\setminus\mathcal{S}$, let $\beta_i$ denote the number of grounded nodes that are in the neighborhood of $v_i$.  We can then write the grounded Laplacian matrix $L_g$ as $L_g=\bar{L}+E$,
where $E = \diag(\beta_1, \beta_2, \ldots, \beta_{n-|\mathcal{S}|})$ and $\bar{L}$ is the Laplacian matrix for the graph induced by the nodes $\mathcal{V}\setminus\mathcal{S}$.  Using Weyl's inequality for $i=1,2,..., n-|\mathcal{S}|$, we have 
\begin{equation}
\lambda_i(\bar{L})\leq \lambda_i(L_g)\leq \lambda_i(\bar{L})+|\mathcal{S}|.
\label{eqn:weylll}
\end{equation}
Thus we have,
\begin{multline}
\frac{1}{2}\left(\sum_{i=2}^{n-|\mathcal{S}|}\left(\frac{1}{\lambda_i(\bar{L})+|\mathcal{S}|}\right)+\frac{1}{\lambda_1(L_g)}\right) \leq \frac{1}{2}\trace(L_g^{-1}) 
\leq \frac{1}{2}\left(\sum_{i=2}^{n-|\mathcal{S}|}\frac{1}{\lambda_i(\bar{L})}+\frac{1}{\lambda_1(L_g)}\right).
\label{eqn:per3}
\end{multline}
Noting that $\bar{L}$ is the Laplacian matrix for an Erdos-Renyi random graph on $n-|\mathcal{S}|$ nodes with constant $p$, for $n-|\mathcal{S}|= \Omega(n)$ we have $(1- o(1))(n-|\mathcal{S}|)p \le \lambda_2(\bar{L})\le(1+o(1))(n-|\mathcal{S}|)p$ and $(1 - o(1))(n-|\mathcal{S}|)p \le \lambda_{n-|\mathcal{S}|}(\bar{L}) \le (1+ o(1))(n-|\mathcal{S}|)p$ asymptotically almost surely \cite{Mesbahi}. Thus according to \eqref{eqn:per3}, Theorem~\ref{thm:cohh} and considering the fact that $|\mathcal{S}|=o(\sqrt{n})$, the result is obtained. 
\end{IEEEproof}

\begin{remark}
With a single leader and constant $p$, both $\mathcal{H}_2$ and $\mathcal{H}_{\infty}$ disorders are within $(1\pm o(1))\frac{1}{p}$ asymptotically almost surely (from \eqref{eqn:per4} and \eqref{eqn:per44}).
\end{remark}

Theorems \ref{thm:erdos1} and \ref{thm:cohh} apply to any edge probability $p$ satisfying $p(n) \ge \frac{c\ln{n}}{n}$, for constant $c>1$.  For any $p$ in this range, the second smallest eigenvalue of the graph Laplacian satisfies $\lambda_2({L})=\Theta(np)$ asymptotically almost surely \cite{PiraniSundaramArxiv}. Thus, for this more general class of edge probabilities, we have the following looser bound on $\mathcal{H}_2$ disorder.
\begin{corollary} 
Consider a random  graph $\mathcal{G}(n,p)$ where the edge probability satisfies $p\ge\frac{c\ln n}{n}$ for any $c>1$ and a set of grounded nodes $\mathcal{S} \subset \mathcal{V}$ chosen uniformly at random such that $|\mathcal{S}|=o(\sqrt{np})$. Then for $\mathcal{H}_{2}$ disorder we have
\begin{equation}
 \frac{1}{2}\trace(L_g^{-1})=\Theta(\frac{1}{p}),
\label{eqn:per5}
\end{equation}
asymptotically almost surely.
\end{corollary}

\begin{IEEEproof}
By Theorem \ref{thm:erdos1} we have $(1-o(1))|\mathcal{S}|p \le \lambda_1(L_g) \le (1+o(1))|\mathcal{S}|p$ for the regime of $p$ mentioned in the corollary. 
According to the Cauchy interlacing theorem and \cite{PiraniSundaramArxiv} we have $\beta np \geq 2d_{max}\geq \lambda_n(L) \geq  \lambda_i(L_g)\geq \lambda_i(L)\geq \lambda_2(L) \geq \alpha np$ asymptotically almost surely for some $\alpha, \beta >0$ and $i= 2, ..., n-|\mathcal{S}|$.  Summing the inverse of these eigenvalues to obtain $\trace(L_g^{-1})$ gives the result. 
\end{IEEEproof}

\subsubsection{Delay Threshold in ER Random Graphs}

The following result discusses the value of $\hat{\tau}_{max}$ in ER random graphs.

\begin{theorem}
Consider a random graph $\mathcal{G}(n,p)$ and let $\mathcal{S} \subset \mathcal{V}$ be a set of grounded nodes chosen uniformly at random with $|\mathcal{S}|=o({np})$. Then for constant $p$ we have 
\begin{equation}
(1 - o(1))\frac{\pi}{2np}\leq \hat{\tau}_{max} \leq (1 + o(1))\frac{\pi}{2np},
\end{equation}
asymptotically almost surely. Moreover, for $p \geq \frac{c\ln n}{n}$ and $c>1$ we have 
\begin{equation}
\hat{\tau}_{max} = \Theta(\frac{1}{np}),
\end{equation}
asymptotically almost surely.
\label{thm:jsbvh}
\end{theorem}

\begin{IEEEproof}
Similar to the proof of Theorem \ref{thm:coh}, we  write the grounded Laplacian matrix $L_g$ as $L_g=\bar{L}+E$.  Using Weyl's inequality \eqref{eqn:weylll} for constant $p$ we have $(1 - o(1))(n-|\mathcal{S}|)p \le \lambda_{n-|\mathcal{S}|}(\bar{L}) \le (1+ o(1))(n-|\mathcal{S}|)p$ asymptotically almost surely \cite{Mesbahi}. By considering $|\mathcal{S}|=o({np})$ we have 
\begin{equation}
(1 - o(1))np \le \lambda_{n-|\mathcal{S}|}({L}_g) \le (1+ o(1))np,
\label{eqn:ght}
\end{equation}
asymptotically almost surely, which yields the result. 

For $p \geq \frac{c\ln n}{n}$ according to Theorem \ref{thm:bounddmax} and considering the fact that in this regime of $p$ we have $d_i=\Theta(np)$ for all $v_i\in \mathcal{V}$ \cite{PiraniSundaramArxiv},  we have $\lambda_{n-|\mathcal{S}|}(\bar{L})=\Theta(np)$  asymptotically almost surely. Moreover, based on Weyl's inequality \eqref{eqn:weylll} and according to the fact that $|\mathcal{S}|=o({np})$ we have $\lambda_{n-|{\mathcal{S}}|}(L_g)=\Theta(np)$, which yields the result. 
\end{IEEEproof}

\subsection{Random Regular Graphs}

\begin{definition}
Let $\Omega_{n,d}$ be the set of all undirected graphs on $n$ nodes where every node has degree $d$ (note that this assumes that $nd$ is even).  A {\it random $d$-regular graph} (RRG), denoted $\mathcal{G}_{n,d}$ is a graph drawn uniformly at random from $\Omega_{n,d}$.  
\end{definition}

\subsubsection{Network disorder in RRG}

We have the following result for the disorder in random regular graphs. 

\begin{theorem} 
Let $\mathcal{G}_{n,d}$ be a random $d$-regular graph on $n$ nodes, with a set of leaders $\mathcal{S}$. Then for sufficiently large (constant) $d$, both $\mathcal{H}_2$ and $\mathcal{H}_{\infty}$  disorders are $O(n)$
asymptotically almost surely.
\label{thm:ry}
\end{theorem}

\begin{IEEEproof}
Let $v_i\in \mathcal{S}$ be any arbitrary vertex in the leader set. It was shown in \cite{PiraniSundaramArxiv} that  $\lambda_1(L_{gi})=\Theta(\frac{1}{n})$ asymptotically almost surely for a random $d$-regular graph with sufficiently large $d$, where $L_{gi}$ is the grounded Laplacian  induced by node $v_i$. Based on the interlacing theorem we have $\lambda_1(L_{g}(\mathcal{S})) \geq \lambda_1(L_{gi})$. This implies that the $\mathcal{H}_{\infty}$ disorder (given by \eqref{eqn:hinfti}) in this case is $O(n)$. Moreover, for  $j=2,3,...,n-1$, we have 
\begin{equation}
 \lambda_j(L_{g}(\mathcal{S})) \geq \lambda_j(L_{gi})\geq \lambda_2(L_{gi}) \ge \lambda_2(L)\geq \alpha d,
 \label{eqn:bunchofinq}
\end{equation}
asymptotically almost surely for some $\alpha>0$ and sufficiently large $d$. The first (left) and the second last inequalities in \eqref{eqn:bunchofinq} are due to the interlacing theorem, and the last inequality is a direct consequence of the result in \cite{Friedman}. Thus for $\mathcal{H}_2$ disorder we have
$\frac{1}{2}\trace(L_{gi}^{-1}) = \frac{1}{2}\sum_{j = 1}^{n-1}\frac{1}{\lambda_j(L_{gi})} = \frac{1}{2\lambda_1(L_{gi})} + \frac{1}{2}\sum_{j = 2}^{n-1}\frac{1}{\lambda_j(L_{gi})} 
= \Theta(n) + O(\frac{n}{d}) = \Theta(n)$. Hence, from the first inequality in \eqref{eqn:bunchofinq} we have $\frac{1}{2}\sum_{j = 1}^{n-1}\frac{1}{\lambda_j(L_{g}(\mathcal{S}))}=O(n)$, which gives the result.
\end{IEEEproof}

Based on Theorem~\ref{thm:ry}, for any combination of  leaders in a random regular graph, both $\mathcal{H}_2$ and $\mathcal{H}_{\infty}$  disorders will grow at most linearly with the network size.\footnote{There are some graphs in which the network disorder grows faster than the network size, e.g., $d$-dimensional grids for $d=1,2$ in which the network disorder is $O(n^2)$ and $O(n\log (n))$, respectively \cite{Bamieh2}.}

\subsubsection{Delay Threshold in RRG}

The following result discusses the value of $\hat{\tau}_{max}$ in random regular  graphs. It follows immediately from Theorem \ref{thm:bounddmax} and \eqref{eqn:ti}  and thus we skip the proof.  

 \begin{theorem} 
Let $\mathcal{G}_{n,d}$ be a random $d$-regular graph on $n$ nodes, with a set of leaders $\mathcal{S}$. Then we have $\frac{\pi}{4d}\leq \hat{\tau}_{max}\leq \frac{\pi}{2d}$.
\end{theorem}

The values of system $\mathcal{H}_2$ and  $\mathcal{H}_{\infty}$ disorders and $\hat{\tau}_{max}$ in ER random graphs and RRGs are summarized in Table \ref{tab:nscvh}.

\begin{table}[h] 
\begin{tabular}{c c c c c} 

\textbf{ER}   \\
 \hline
 $p$ & $\mathcal{H}_2$ &  $\mathcal{H}_{\infty}$ & $\hat{\tau}_{max}$\\
 \hline
 Constant  & $(1\pm o(1))\frac{|\mathcal{S}|+1}{2|\mathcal{S}|p}$ & $(1\pm o(1))\frac{1}{|\mathcal{S}|p}$ & $(1\pm o(1))\frac{\pi}{2np}$\\
 $\frac{c\ln n}{n},c>1$ & $\Theta(\frac{1}{p})$ & $(1\pm o(1))\frac{1}{|\mathcal{S}|p}$ & $\Theta(\frac{1}{np})$\\
 \hline\\ \\
\textbf{RRG}  \\ 
\hline
$d$ & $\mathcal{H}_2$ &  $\mathcal{H}_{\infty}$ & $\hat{\tau}_{max}$\\
 \hline
Constant & $O(n)$ & $O(n)$ & $\in [\frac{\pi}{4d},\frac{\pi}{2d}]$\\
\end{tabular}
\centering
\caption {Disorder and Delay Threshold in random graphs. For network disorders in ER graphs, it is assumed that $|\mathcal{S}|=o(\sqrt{np})$.}
\label{tab:nscvh}
\end{table}

\section{System Robustness as Network Centrality Metrics} 
\label{sec:robustleader}

In this section, we look at the system robustness to disturbances and time delay via network centrality metrics. In particular, we seek to choose a leader  in order to maximize robustness to disturbances \cite{Bamieh}, \cite{Bamieh2}, \cite{Leonard}, \cite{piranicdc} and time delay.  As argued in Section~\ref{sec:robust_disturbances}, minimizing the $\mathcal{H}_{\infty}$ disorder is equivalent to maximizing convergence rate, and thus we will make connections to existing work that has looked at this latter metric \cite{Clark1, Ghaderi13,piranicdc}.

\subsection{Optimal Leaders for Each Objective}

In this subsection, we provide conditions for a single leader in a network to optimize each robustness metric separately. 

\subsubsection{Minimizing $\mathcal{H}_2$ disorder}

As shown in \cite{Leonard}, the  optimal single leader for minimizing  $\mathcal{H}_2$ disorder  is the node  with maximal {\it information centrality} defined as $IC(\mathcal{G})=\max_{i\in \mathcal{V}}[\frac{1}{n}\sum_j\gamma_{ij}]^{-1}$,
where $\gamma_{ij}$ is the sum of the lengths of {\it all} paths  between nodes $v_i$ and $v_j$ in the network.
 For the case of trees, the information central vertex and the closeness central vertex (a vertex whose summation of distances to the rest of the vertices is minimum)  in the network are the same. The result is extended to the case of multiple leaders in \cite{Fitch2}.

\subsubsection{Minimizing $\mathcal{H}_{\infty}$ disorder}

In order to discuss optimal leaders for minimizing the $\mathcal{H}_{\infty}$ disorder (or maximizing the convergence rate), we define a {\it grounding centrality} for each node $v_i$, which is equal to the smallest eigenvalue of the grounded Laplacian induced by that node.  Thus, the single best leader in terms of minimizing $\mathcal{H}_{\infty}$ disorder (maximizing convergence rate) is the one with largest grounding centrality \cite{ACC}. 

In the following example, we show that the grounding central vertex and the information central vertex can be far from each other in a graph and consequently the best leader to minimize $\mathcal{H}_{2}$ disorder can be different from the one that minimizes $\mathcal{H}_{\infty}$ disorder.

\begin{example}
A {\it broom tree}, $B_{n,\Delta}$, is a star  $S_{\Delta}$  with $\Delta$ leaf vertices and a path of length $n-\Delta-1$ attached to the center of the star, as illustrated in Fig.~\ref{fig:broom}  \cite{Stevanovic10}.
\begin{figure}[hb]
\begin{center}
\begin{tikzpicture}
[scale=.20, inner sep=1.5pt, minimum size=3pt,auto=center,every node/.style={circle, draw=black, thick}]
\node (n1) at (7,7) {1};
\node (n2) at (1,7) {2};
\node (n3) at (1,1) {3};
\node (n4) at (7,1) {4};
\node (n5) at (4,4) {5};
\node (n6) at (9,4) {6};
\node (n7) at (14,4) {7};
\node (n8) at (19,4) {8};
\node (n9) at (24,4) {9};
\foreach \from/\to in {n1/n5,n5/n2,n3/n5,n5/n4,n6/n5,n7/n6,n7/n8,n8/n9}
\draw (\from) -- (\to);
\end{tikzpicture}
\end{center}
\caption{Broom tree with $\Delta=4$, $n=9$.}
\label{fig:broom}
\end{figure}
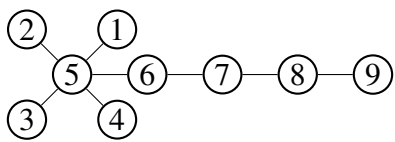

Consider the broom tree $B_{2\Delta+1,\Delta}$.  By numbering the vertices as shown in Fig.~\ref{fig:broom}, for $\Delta=500$, we find (numerically) that the grounding central vertex is vertex 614.  The information central vertex is located at the middle of the star (vertex 501). The deviation of the grounding central vertex from the information central vertex increases as $\Delta$ increases.\footnote{In this example, some other well known centrality metrics, e.g. degree centrality, betweenness centrality, closeness centrality and eigenvector centrality,  are all optimized at the center of the star. }
\label{exm:1}
\end{example}

\subsubsection{Maximizing Delay Threshold $\hat{\tau}_{max}$}

In the following lemma, we provide a sufficient condition for a leader in a network to maximize the delay threshold $\hat{\tau}_{max}$. 

\begin{lemma}
Consider a connected graph $\mathcal{G} = \{\mathcal{V},\mathcal{E}\}$.  Node $v_k \in \mathcal{V}$ is the optimal leader for maximizing $\hat{\tau}_{\max}$ if $d_k\geq 2d_i$ for all $v_i\in \mathcal{V}\setminus \{v_k\}$.
\label{lem:aergtn}
\end{lemma}

\begin{IEEEproof}
 Based on \eqref{eqn:ti}, we need to  show that 
\begin{equation}
\lambda_{n-|\mathcal{S}|}(L_{gk})\leq \lambda_{n-|\mathcal{S}|}(L_{gi}),
\label{eqn:sgb}
\end{equation}
 for all $i\in \mathcal{V}\setminus \{v_k\}$. Here $L_{gi}$ and $L_{gk}$ are the grounded Laplacian matrices induced by nodes $v_i$ and $v_k$, respectively. Based on Theorem \ref{thm:bounddmax}, a  sufficient condition for \eqref{eqn:sgb} is  
\begin{equation}
\max_{{u\in \mathcal{V}\setminus\{v_i\}}}d_u \geq \max \left\{\max_{{u\in \mathcal{V}\setminus\{v_k\}}}d_u, \max_{\substack{(u,v)\in \mathcal{E}\\{u,v\in \mathcal{V}\setminus\{v_k\}}}}\{d_u+d_v\}\right\},
\label{eqn:kohig}
\end{equation}
for all $v_i\in \mathcal{V}$. We know that 
$$2 \max_{u\in \mathcal{V}\setminus\{v_k\}}d_u \geq \max_{\substack{(u,v)\in \mathcal{E}\\{u,v\in \mathcal{V}\setminus\{v_k\}}}}\{d_u+d_v\}.$$
Thus, a sufficient condition for \eqref{eqn:kohig} is 
\begin{equation}
\max_{{u\in \mathcal{V}\setminus\{v_i\}}}d_u \geq 2 \max_{u\in \mathcal{V}\setminus\{v_k\}}d_u,
\end{equation}
for all $v_i\in \mathcal{V}\setminus \{v_k\}$, which is equivalent to $d_k\geq 2d_i$ for all $v_i\in \mathcal{V}\setminus \{v_k\}$.
\end{IEEEproof}

Based on Lemma \ref{lem:aergtn},  the leader which maximizes $\hat{\tau}_{max}$ in Example \ref{exm:1} is the center of the star, which is the same leader which minimizes $\mathcal{H}_2$ disorder. However it is not true that these two robustness metrics ($\mathcal{H}_2$ disorder and $\hat{\tau}_{max}$) always share the same optimal leader. For example, in the graph shown in Example \ref{exm:1}, if we fix the degree of the star and increase the length of the tail, the information central vertex will no longer remain in the center of the star, while the center of the star still maximizes $\hat{\tau}_{max}$, provided its degree is at least 4.

Lemma \ref{lem:aergtn} indicates that if a node in a network has substantially higher degree than the other nodes, then it is  the optimal leader for $\hat{\tau}_{max}$. However the fact that the highest degree node is always the best leader is not true when the differences in degrees are moderate, as shown in the following example.

\begin{example}
In the graph shown in  Fig.~\ref{fig:brooms}, the black nodes have degree 3 which is the highest degree in the graph. However,  we have $\lambda_{n-|\mathcal{S}|}^{Gray}=3.7321$ and $\lambda_{n-|\mathcal{S}|}^{Black}=4.1149$, where $\lambda_{n-|\mathcal{S}|}^{Gray}$ is the largest eigenvalue of the grounded Laplacian induced by the gray node (and the same for $\lambda_{n-|\mathcal{S}|}^{Black}$).

\begin{figure}[hb]
\begin{center}
\begin{tikzpicture}
[scale=.15, inner sep=1.5pt, minimum size=8pt,auto=center,every node/.style={circle, draw=black, thick}]

\node (n1) at (11,7) {};
\node (n2) at (11,1) {};
\node (n3) at (15,4) {};
\node (n4) at (21,4) {};
\node (n5) at (27,4) {};
\node (n6) at (31,7) {};
\node (n7) at (31,1) {};
\foreach \from/\to in {n1/n3, n3/n2, n3/n4, n5/n4, n5/n6, n5/n7}
\draw (\from) -- (\to);
\filldraw[color=black!100, fill=black!100, thick](15,4) circle (0.8);
\filldraw[color=black!100, fill=black!100, thick](27,4) circle (0.8);
\filldraw[color=black!100, fill=gray!100, thick](21,4) circle (0.8);

\end{tikzpicture}
\end{center}
\caption{An example which shows that a leader with maximum degree does not maximize $\hat{\tau}_{max}$.}
\label{fig:brooms}
\end{figure}
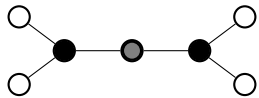

\label{exm:2}
\end{example}

Our discussion and results in this section have shown that the optimal leaders for each of the robustness metrics will, in general, be different. In the following subsection, we discuss  conditions under which a single leader can optimize all of the three objectives ($\mathcal{H}_{\infty}$ disorder or convergence rate, $\mathcal{H}_{2}$ disorder, and delay threshold) simultaneously.

\subsection{A Sufficient Condition for a Leader to Minimize Network Disorder and Maximize Delay Threshold}

In this subsection we provide a sufficient condition for a leader to simultaneously minimize $\mathcal{H}_2$ and $\mathcal{H}_{\infty}$ disorders and maximize $\hat{\tau}_{max}$.  We require the following concept.

\begin{definition}
 The {\it resistance distance} $r_{ij}$ between two vertices $v_i$ and $v_j$ in a graph is the equivalent resistance between these two vertices when we treat each edge of the graph as a $1 \Omega$ resistor. The {\it effective resistance} of vertex $v_i$ is $R_i=\sum_{j\neq i}r_{ij}$. \end{definition}

 It can be shown that the resistance distance between $v_i$ and $v_j$ is the $j$-th diagonal element of $L_{gi}^{-1}$, where $v_i$ is a single grounded vertex \cite{Bapat}.  Thus, the effective resistance of vertex $v_i$ is 
\begin{equation}
R_i=\trace(L_{gi}^{-1}).
\label{eqn:dop1}
\end{equation}
Moreover,  the resistance distance between vertices $v_i$ and $v_j$  is given by \cite{Bapat}
\begin{equation}
r_{ij}=(e_i-e_j)^T L_{gk}^{-1} (e_i-e_j),
\label{eqn:resis}
\end{equation}
where $k \notin \{i,j\}$ is the index of an arbitrary vertex  which becomes grounded and $e_i$ is a vector of zeros except for a $1$ in the element corresponding to the  $i$-th vertex.

\begin{theorem}
Consider a connected graph $\mathcal{G} = \{\mathcal{V},\mathcal{E}\}$.  Node $v_k \in \mathcal{V}$ will simultaneously be the best single leader to minimize $\mathcal{H}_{\infty}$ disorder, minimize $\mathcal{H}_2$ disorder, and maximize $\hat{\tau}_{max}$ and convergence rate if $d_k \ge \frac{2d_i}{x^2_{min}}$ for all $v_i \in \mathcal{V}\setminus\{v_k\}$, where $x_{min}$ is the smallest component of a nonnegative eigenvector corresponding to the smallest eigenvalue of the grounded Laplacian $L_{gk}$.
\label{thm:jvfvv}
\end{theorem}

\begin{IEEEproof}
From \eqref{eqn:dop1} and \eqref{eqn:resis}, the effective resistance of $v_i$  is
  \begin{align}
  R_i = \trace(L_{gi}^{-1}) &=r_{ik}+\sum_{j\neq i}(e_i-e_j)^T L_{gk}^{-1} (e_i-e_j)\nonumber \\
   &= \trace(L_{gk}^{-1})+nr_{ik}-2S_{i}^k,
  \label{eqn:resis2}
  \end{align}
  where $S_{i}^k$ is the sum of the elements of the $i$-th row (or column) in $L_{gk}^{-1}$. From \eqref{eqn:resis2} we have $\trace(L_{gk}^{-1})-\trace(L_{gi}^{-1})=2S_{i}^k-nr_{ik}$. Thus for $v_k$ to be a better leader than $v_i$ for network $\mathcal{H}_2$ disorder, it is sufficient to have
  \begin{equation}
  2\bar{S}-nr_{ik}\leq 0,
  \label{eqn:resis4}
  \end{equation}
  where $\bar{S}=\max_j S_{j}^k$ is the maximum row sum in $L_{gk}^{-1}$. On the other hand, from \cite{ACC} we know that $\bar{S}x_{min}\leq \lambda_{max}(L_{gk}^{-1})\leq \bar{S}$. Combining this with \eqref{eqn:resis4} yields $\lambda_{max}(L_{gk}^{-1})\leq \frac{nr_{ik}x_{min}}{2}$
as a sufficient condition for $v_k$ to be a better leadership candidate than $v_i$ for the objective of minimizing $\mathcal{H}_2$ network disorder. This sufficient condition can be more conveniently framed as $\lambda_1(L_{gk}) \geq \frac{2}{nr_{ik}x_{min}}$.
  From \cite{Bapat} we know that $r_{ik}\geq \max\{\frac{1}{d_i},\frac{1}{d_k}\}$ where $d_i$ and $d_k$ are the degrees of vertices of $v_i$ and $v_k$ respectively. Thus a sufficient condition for the above inequality to hold is
$\lambda_1(L_{gk}) \geq \frac{2\min\{d_i,d_k\}}{nx_{min}}$.
A sufficient condition for this, based on \eqref{eqn:maineq} (with $\mathcal{S}=\{v_k\}$), is
\begin{equation}
\frac{d_kx_{min}}{n-1} \geq \frac{2\min\{d_i,d_k\}}{nx_{min}}.
\label{eqn:obj}
\end{equation}
On the other hand, for $v_k$ to be a better leader compared to $v_i$ for optimizing $\mathcal{H}_{\infty}$ disorder (or equivalently maximizing convergence rate),   according to \eqref{eqn:maineq} it is sufficient to have $\frac{d_k x_{min}}{n-1} \geq \frac{d_i}{n-1}$
which gives $d_k \geq \frac{d_i}{x_{min}}$, where $x_{min}$ is again the smallest eigenvector component of $\mathbf{x}$, a nonnegative eigenvector corresponding to $ \lambda_1(L_{gk})$. Combining this with \eqref{eqn:obj}, a sufficient condition for $v_k$ to be a better leader than $v_i$ for both objectives simultaneously is 
\begin{equation}
d_k \geq \max\left\{\frac{d_i}{x_{min}}, \frac{2d_i(n-1)}{nx^2_{min}}\right\}.
\end{equation}
Since $\frac{2(n-1)}{nx_{min}} \ge 1$ for $n \ge 2$, we conclude that $d_k \ge \frac{2d_i}{x^2_{min}}$ for all $v_i \in \mathcal{V}\setminus\{v_k\}$ is a sufficient condition for $v_k$ to be the optimal leader for minimizing $\mathcal{H}_2$ and $\mathcal{H}_{\infty}$ disorders simultaneously. 

Based on Lemma \ref{lem:aergtn} and considering the fact that  $\frac{2d_i}{x^2_{min}}\geq 2d_i$, the condition mentioned in the theorem is also a sufficient condition for $d_k\geq 2d_i$ and thus it is sufficient for $v_k$ to be a leader which maximizes $\hat{\tau}_{max}$, which yields the result. 
\end{IEEEproof}

\begin{remark}
Based on the sufficient condition $d_k \ge \frac{2d_i}{x^2_{min}}$ and the bound given in Lemma \ref{lem:mainlem} for $x_{min}$, as the algebraic connectivity of the graph induced by the follower agents becomes larger, $x_{min}$ becomes closer to 1 and  the condition $d_k \ge \frac{2d_i}{x^2_{min}}$ will be less demanding in terms of the degree that the leader agent $d_k$ is required to have. However, the condition for the optimal leader for $\hat{\tau}_{max}$ given in Lemma \ref{lem:aergtn} is independent of the connectivity of the graph induced by the follower agents. 
\end{remark}

In the following example, we describe a graph such that  there exists a node satisfying $d_k \ge \frac{2d_i}{x^2_{min}}$ and consequently becomes the optimal leader to optimize all of the objectives, i.e.  $\mathcal{H}_{2}$ disorder, $\mathcal{H}_{\infty}$ disorder (and convergence rate), and delay threshold $\hat{\tau}_{max}$.

\begin{example}
Consider an ER random graph $\mathcal{G}(n,p)$ with $p \geq \frac{c\ln n}{n}$ for some $c>1$. The degree of each vertex in the graph is $d_i=\Theta(np)$  and the algebraic connectivity is $\lambda_2(L) = \Theta(np)$ asymptotically almost surely \cite{PiraniSundaramArxiv}.  Suppose we wish to connect a single leader node $v_{n+1}$ to this network in such a way that it is the single best leader for optimizing $\mathcal{H}_2$ disorder, $\mathcal{H}_{\infty}$ disorder (or convergence rate) and the delay threshold.   Pick any $\epsilon > 0$ and connect $v_{n+1}$ to any $(2+\epsilon)d_{max}$ nodes in the network, where $d_{max}$ is the maximum degree of any node in the network.  Let $L_{g,n+1}$ be the grounded Laplacian induced by $v_{n+1}$. From \eqref{eqn:eige12}, the eigenvector for $\lambda_1(L_{g,n+1})$ has smallest component $x_{min} \ge 1 - \frac{2\sqrt{d_{n+1}}}{\lambda_2({L})}$ which goes to $1$ asymptotically almost surely (since $d_{n+1} = (2+\epsilon) d_{max} = \Theta(np)$ and $\lambda_2(L) = \Theta(np)$).  Thus the condition $d_{n+1} \ge \frac{2d_i(n-1)}{nx^2_{min}}$ will be satisfied for this node asymptotically almost surely.  
\end{example}

\section{Simulation Results}
\label{sec:simulations}

\begin{figure}[h!]
\centering
\includegraphics [width=8cm, height=5cm]{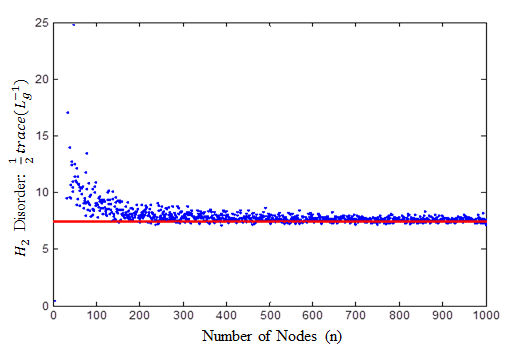}
\caption{$\mathcal{H}_2$ disorder in $\mathcal{G}(n,0.1)$ converges to 7.5, the red line.}
\label{fig:h2dis}
\end{figure}

\begin{figure}[h!]
\centering
\includegraphics[width=8cm, height=5cm]{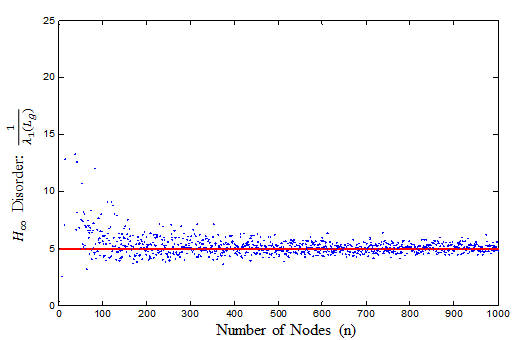}
\caption{$\mathcal{H}_{\infty}$ disorder in $\mathcal{G}(n,0.1)$ converges to 5, the red line.}
\label{fig:hinfdis}
\end{figure}

\begin{figure}[h!]
\centering
\includegraphics[width=8.2cm, height=5cm]{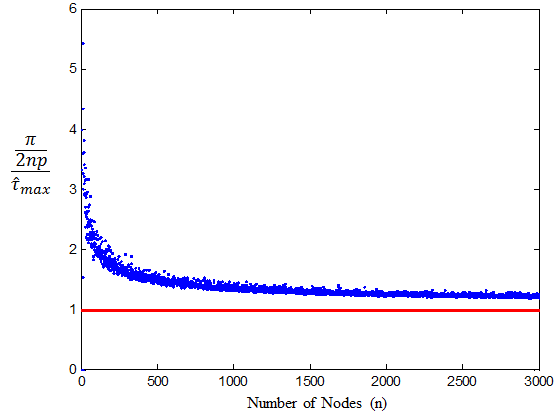}
\caption{$\frac{\frac{\pi}{2np}}{\hat{\tau}_{max}}$  in $\mathcal{G}(n,0.1)$ converges to $1$, shown in red.}
\label{fig:delayfig}
\end{figure}

In this section, we provide some simulation results for the robustness of the leader-follower consensus dynamics when the underlying network is an  ER random graph $\mathcal{G}(n,p)$ for constant $p$. We set the number of leaders in a network as $|\mathcal{S}|=2$ and the probability of edge formation as $p=0.1$. In Fig. \ref{fig:h2dis}, the value of the $\mathcal{H}_2$ disorder converges to $\frac{|\mathcal{S}|+1}{2|\mathcal{S}|p}=7.5$, as predicted by Theorem \ref{thm:coh}. In Fig. \ref{fig:hinfdis}, the value of the $\mathcal{H}_{\infty}$ disorder converges to $\frac{1}{|\mathcal{S}|p}=5$, in accordance with  Theorem \ref{thm:cohh}.  In Fig. \ref{fig:delayfig},  $\frac{\frac{\pi}{2np}}{\hat{\tau}_{max}}$  converges to $1$ (although slowly),  as specified by Theorem \ref{thm:jsbvh}.

